\documentclass[11pt]{article}
\usepackage{enumerate}
\usepackage{amssymb,a4wide,latexsym,makeidx,epsfig,fleqn}
\usepackage{amsthm}
\usepackage{amsmath}
\usepackage{enumerate}
\usepackage{mathrsfs}

\newtheorem{theorem}{Theorem}[section]

\newtheorem{lemma}[theorem]{Lemma}
\newtheorem{fact}[theorem]{Fact}

\newtheorem{proposition}[theorem]{Proposition}

\begin{document}
\textwidth 150mm \textheight 225mm
\title{Gallai-Ramsey numbers for a class of graphs with five vertices
\thanks{Supported by the National Natural Science Foundation of China (No. 11871398) and
the Natural Science Basic Research Plan in Shaanxi Province of China (Program No. 2018JM1032).}}
\author{{Xihe Li$^{a,b}$ and Ligong Wang$^{a,b,}$\thanks{Corresponding
author.}}\\
{\small $^a$Department of Applied Mathematics, School of Science,}\\ {\small Northwestern
Polytechnical University, Xi'an, Shaanxi 710072,
P. R. China. }\\ {\small $^b$Xi'an-Budapest Joint Research Center for Combinatorics,}\\ {\small Northwestern Polytechnical
University, Xi'an, Shaanxi 710129, P. R. China.
}\\ {\small E-mail: lxhdhr@163.com; lgwangmath@163.com}}
\date{}
\maketitle
\begin{center}
\begin{minipage}{120mm}
\vskip 0.3cm
\begin{center}
{\small {\bf Abstract}}
\end{center}
{\small Given two graphs $G$ and $H$, the $k$-colored Gallai-Ramsey number $gr_k(G : H)$
is defined to be the minimum integer $n$ such that every $k$-coloring of the
complete graph on $n$ vertices contains either a rainbow copy of $G$ or a monochromatic copy
of $H$. In this paper, we consider $gr_k(K_3 : H)$ where $H$ is a connected graph with five vertices
and at most six edges. There are in total thirteen graphs in this graph class, and the
Gallai-Ramsey numbers for some of them have been studied step by step in several papers.
We determine all the Gallai-Ramsey numbers for the remaining graphs, and we also obtain some
related results for a class of unicyclic graphs.

 }
\end{minipage}
\end{center}

\section{Introduction}
\label{sec:ch-introduction}

In this paper, we only consider edge-colorings of finite simple graphs.
For an integer $k \geq 1$, let $c$ :
$E(G)\rightarrow [k]$ be a $k$-edge-coloring (or simply, a $k$-coloring)
of a graph $G$, where $[k] := \{1, 2, \ldots, k\}$.
Note that $c$ is not necessarily a proper edge-coloring.
A coloring of a graph is called {\it monochromatic} if all edges are colored the same,
and a coloring is called {\it rainbow} if all edges are colored differently.
A graph $G$ is called {\it rainbow triangle-free} if every triangle in $G$ contains at most two colors.
A {\it Gallai-$k$-coloring} is a rainbow triangle-free coloring of a complete graph using at
most $k$ colors. In 1967, Gallai \cite{Gallai} provided the following important
structural result in his original paper, which was translated
into English and endowed by comments in \cite{MaPr}.

\noindent\begin{theorem}\label{th:Gallai} {\normalfont (\cite{Gallai,MaPr})}
In any rainbow triangle-free coloring of a complete graph, there exists
a partition $V_1$, $V_2$, $\ldots$ , $V_m$ ($m\geq 2$) of the vertices such that between the parts there are in total
at most two colors, and between every pair of parts, there is only one color on the edges.
\end{theorem}

The partition provided in Theorem \ref{th:Gallai} is called a {\it Gallai-partition}.
Given a Gallai-coloring of a complete graph and a Gallai-partition $V_1$, $V_2$, $\ldots$ , $V_m$,
the subgraph induced by $\{v_1, v_2, \ldots, v_m\}$ is called the {\it reduced graph}
of that partition, where $v_i \in V_i$ for every $i \in [m]$.

Given a graph $G$, the {\it Ramsey number} $r_{k}(G)$ is the minimum integer $n$ such that
every $k$-coloring of $K_n$ contains a monochromatic copy of $G$.
Given two graphs $G$ and $H$, the {\it $k$-colored Gallai-Ramsey number} $gr_k(G : H)$
is defined to be the minimum integer $n$ such that every $k$-coloring of the complete graph on $n$
vertices contains either a rainbow copy of $G$ or a monochromatic copy of $H$.
It is clear that $gr_k(G : H) \leq r_{k}(H)$ for any two graphs $G$ and $H$.

Note that finding the exact values of $gr_k(G : H)$ is far from trivial, even for
some graphs $H$ with small number of vertices and edges.
The Gallai-Ramsey numbers $gr_k(K_3 : H)$ for all the graphs $H$ with
four vertices have been found in several previous papers. However, this number for some graphs with
five vertices have not been determined, and finding the exact values for these graphs
is a fundamental work in this research area.
In this paper, we mainly focus on $gr_k(K_3 : H)$ where $H$ is a connected graph with five vertices
and at most six edges. There are thirteen graphs in this graph class (see Fig. 1.1),
and the Gallai-Ramsey numbers for $F_1, F_2, \ldots, F_8$ have been studied step by step in several papers
\cite {FGJM,FuMa1,GSSS,Mag,WMMZ,WMNX}. We determine the Gallai-Ramsey numbers for all the remaining graphs.
For more information on this topic,
we refer the readers to \cite{ChLP,ChGr,FuMa2,GySi,HMOT,SoZh} and two excellent surveys \cite{FuMO1,FuMO2}.

\begin{figure}[htb]
\begin{center}
\includegraphics [width=0.9\textwidth, height=0.4\textwidth]{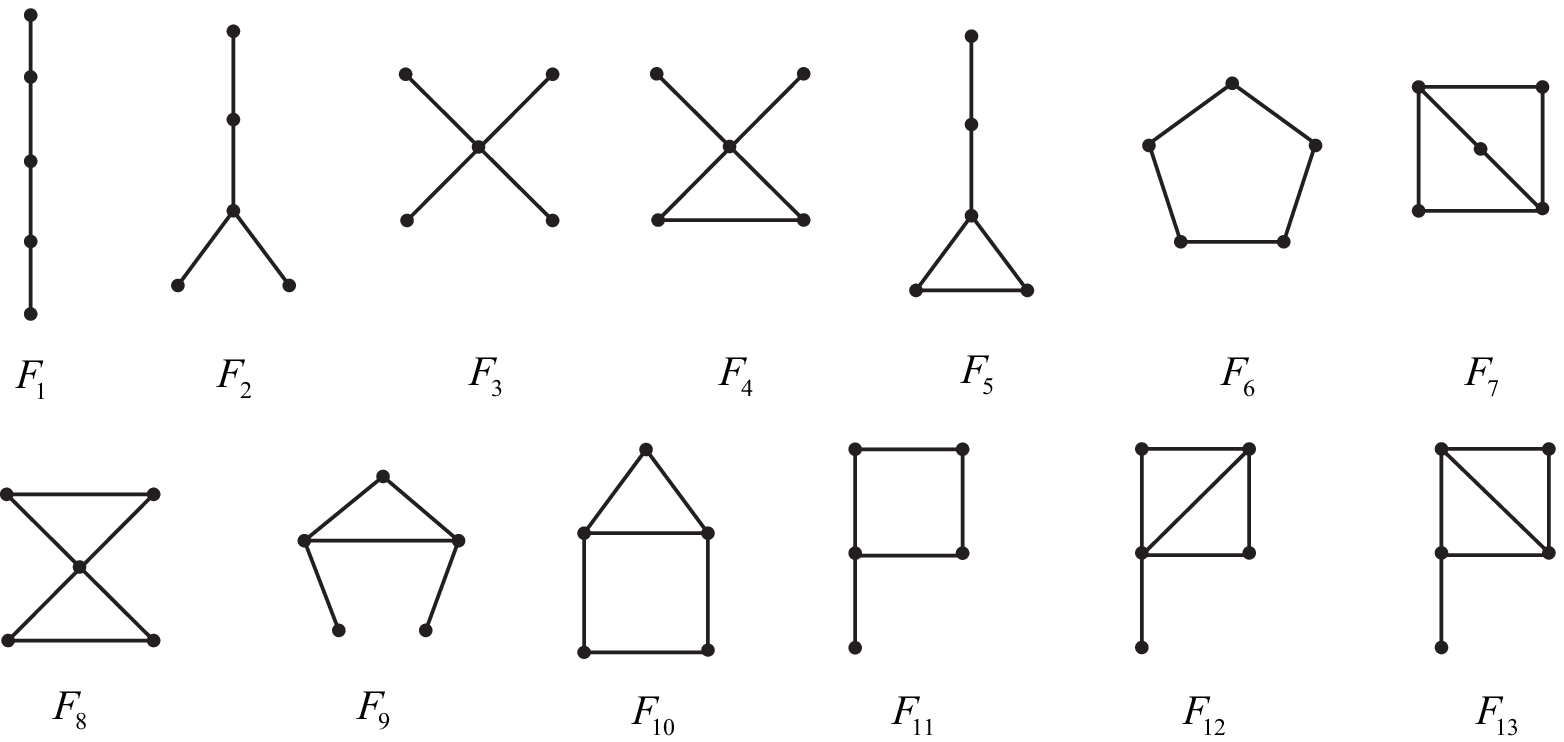}
\centerline{Fig. 1.1: The graphs with five vertices and at most six edges.}
\end{center}
\end{figure}

For graphs $F_9$ and $F_{10}$, we prove the following result in Section \ref{sec:ch-2}.

\noindent\begin{theorem}\label{th:main1} For any integer $k\geq 1$,
$$gr_k(K_3 : F_9)=gr_k(K_3 : F_{10})=
\left\{
   \begin{aligned}
    &8\cdot 5^{(k-2)/2}+1,\ \ \  \ \ \ \ \ & \ & \ \ \mbox{if $k$ is even},\\
    &4\cdot 5^{(k-1)/2}+1,\ \ \ & \ & \ \ \mbox{if $k$ is odd}.
   \end{aligned}
   \right.$$
\end{theorem}

For graph $F_{11}$, we prove the following result in a more general form in Section \ref{sec:ch-3}.
In fact, we consider the Ramsey numbers and Gallai-Ramsey numbers for a class of unicyclic graph $F_{2,n}$.

\noindent\begin{theorem}\label{th:main2} For any integer $k\geq 1$, $gr_k(K_3 : F_{11})=k+4.$
\end{theorem}

For graphs $F_{12}$ and $F_{13}$, we prove the following result in Section \ref{sec:ch-4}.

\noindent\begin{theorem}\label{th:main3} For any integer $k\geq 1$,
$$gr_k(K_3 : F_{12})=gr_k(K_3 : F_{13})=
\left\{
   \begin{aligned}
    &9\cdot 5^{(k-2)/2}+1,\ \ \  \ \ \ \ \ & \ & \ \ \mbox{if $k$ is even},\\
    &4\cdot 5^{(k-1)/2}+1,\ \ \ & \ & \ \ \mbox{if $k$ is odd}.
   \end{aligned}
   \right.$$
\end{theorem}

Furthermore, we define some notation and terminology. Let $c(uv)$ denote
the color used on the edge $uv$. For $U$, $V\subseteq V(G)$ with
$U\cap V= \emptyset$, let $E(U, V)$ denote the set of edges between $U$ and $V$,
and $C(U, V)$ denote the set of colors used on the edges in $E(U, V)$.
If $|C(U, V)|=1$, i.e., all the edges in $E(U, V)$ have a single color, then we use
$c(U, V)$ to denote this color.
In the special case when $U=\{u\}$, we simply write $E(u, V)$,
$C(u, V)$ and $c(u, V)$ for $E(\{u\}, V)$, $C(\{u\}, V)$ and $c(\{u\}, V)$, respectively.
For $v \in V(G)$, let $G-v$ denote the graph obtained from $G$
by deleting the vertex $v$ and the edges incident to $v$.
For $U\subseteq V(G)$, let $G-U$ denote the graph obtained from $G$
by deleting all the vertices of $U$ and the edges incident to some vertex of $U$.
Moreover, let $G[U]$ denote the subgraph of $G$ induced by $U$,
and $C(U)$ denote the set of colors used on the edges of $G[U]$.
For any $i \in C(U)$, the subgraph induced by color $i$ means an induced subgraph
of $G[U]$ on the edges using color $i$.

\section{Proof of Theorem \ref{th:main1}}
\label{sec:ch-2}

We begin with the following simple results which will be used later.

\noindent\begin{proposition}\label{prop:2.1}
If $H_1$ is a subgraph of $H_2$, then $gr_k(K_3 : H_1) \leq gr_k(K_3 : H_2)$.
\end{proposition}

\noindent\begin{lemma}\label{le:2.1} {\normalfont (\cite{Hen})}
$r_2(F_9)=r_2(F_{10})=9$.
\end{lemma}

\noindent\begin{lemma}\label{le:K3} {\normalfont (\cite{ChGr})} For any integer $k\geq 1$,
$$gr_k(K_3 : K_3)=
\left\{
   \begin{aligned}
    &5^{k/2}+1,\ \ \  \ \ \ \ \ & \ & \ \ \mbox{if $k$ is even},\\
    &2\cdot 5^{(k-1)/2}+1,\ \ \ & \ & \ \ \mbox{if $k$ is odd}.
   \end{aligned}
   \right.$$
\end{lemma}

Since $F_9$ is a subgraph of $F_{10}$, we have $gr_k(K_3 : F_9) \leq gr_k(K_3 : F_{10})$ by Proposition \ref{prop:2.1}.
Thus Theorem \ref{th:main1} follows from the following two lemmas immediately.

\noindent\begin{lemma}\label{le:2.2}
For any integer $k\geq 1$,
$gr_k(K_3 : F_9)>
\left\{
   \begin{aligned}
    &8\cdot 5^{(k-2)/2},\ \ \ & \ & \ \ \mbox{if $k$ is even},\\
    &4\cdot 5^{(k-1)/2},\ \ \ & \ & \ \ \mbox{if $k$ is odd}.
   \end{aligned}
   \right.$
\end{lemma}

\noindent {\bf Proof.} If $k$ is even, then let $G_2$ be a 2-colored $K_8$
containing no monochromatic $F_9$ using colors 1 and 2. Suppose that $2i < k$ and we have constructed
a $2i$-coloring $G_{2i}$ of $K_{n_{2i}}$ containing neither rainbow $K_3$ nor monochromatic $F_9$,
where $n_{2i}=8\cdot 5^{(2i-2)/2}$. Let $G'$ be a 2-colored $K_5$ using colors $2i+1$ and $2i+2$
which contains no monochromatic $K_3$, i.e., colors $2i+1$ and $2i+2$ induce two monochromatic
$C_5$. We construct $G_{2i+2}$ by substituting five copies of $G_{2i}$ into vertices of $G'$,
i.e., $G_{2i+2}$ is a blow-up of $G'$. Finally, we obtain a $k$-coloring $G_k$ of $K_n$
containing neither rainbow $K_3$ nor monochromatic $F_9$, where $n=8\cdot 5^{(k-2)/2}$.

If $k$ is odd, then let $G_1$ be a monochromatic $K_4$ using color 1.
Suppose that $2i-1 < k$ and we have constructed
a $(2i-1)$-coloring $G_{2i-1}$ of $K_{n_{2i-1}}$ containing neither rainbow $K_3$ nor monochromatic $F_9$,
where $n_{2i-1}=4\cdot 5^{(2i-2)/2}$. Let $G''$ be a 2-colored $K_5$ using colors $2i$ and $2i+1$
which contains no monochromatic $K_3$. We construct $G_{2i+1}$ by
substituting five copies of $G_{2i-1}$ into vertices of $G''$.
Finally, we obtain a $k$-coloring $G_k$ of $K_n$ containing neither
rainbow $K_3$ nor monochromatic $F_9$, where $n=4\cdot 5^{(k-1)/2}$.
\hfill$\blacksquare$

\noindent\begin{lemma}\label{le:2.3}
For any integer $k\geq 1$,
$gr_k(K_3 : F_{10})\leq
\left\{
   \begin{aligned}
    &8\cdot 5^{(k-2)/2}+1,\ \ \  & \ & \ \ \mbox{if $k$ is even},\\
    &4\cdot 5^{(k-1)/2}+1,\ \ \ & \ & \ \ \mbox{if $k$ is odd}.
   \end{aligned}
   \right.$
\end{lemma}

\noindent {\bf Proof.} We prove the statement by induction on $k$.
The case $k=1$ is trivial and the case $k=2$ holds by Lemma \ref{le:2.1},
so we may assume $k \geq 3$. Suppose $G$ is a Gallai-$k$-coloring of $K_{n}$
containing no monochromatic $F_{10}$, where
$$n=
\left\{
   \begin{aligned}
    &8\cdot 5^{(k-2)/2}+1,\ \ \  & \ & \ \ \mbox{if $k$ is even},\\
    &4\cdot 5^{(k-1)/2}+1,\ \ \ & \ & \ \ \mbox{if $k$ is odd}.
   \end{aligned}
   \right.$$
Let $V_1, V_2, \ldots, V_m$ be a Gallai-partition with $|V_1|\geq |V_2| \geq \cdots \geq |V_m|$
and $m \geq 2$. Since $r_2(F_{10})=9$, we have $2 \leq m \leq 8$. We choose such a partition with $m$ minimum.
Suppose colors 1 and 2 are the two colors used between the parts.
We first have the following simple facts since $G$ is monochromatic $F_{10}$-free.

\noindent\begin{fact}\label{fa:3.1} For any two non-empty parts $V_i$, $V_j$ with $c(V_i, V_j)=\alpha$, we have

(1) there is no monochromatic $P_4$ using color $\alpha$ within $V_i$ and $V_j$, respectively;

(2) if $|V_i|\geq 3$ and $|V_j|\geq 2$, then $\alpha \notin C(V_i)$;

(3) if $|V_i|\geq 2$ and $|V_j|\geq 2$, then there is no vertex with color $\alpha$ to both
$V_i$ and $V_j$ in $V(G)\setminus (V_i \cup V_j)$.
\end{fact}

We now claim that $4 \leq m \leq 8$. In fact, if $2\leq m \leq 3$, then there exists a bipartition
$V_1$ and $V_2$ of $V(G)$ with $|C(V_1, V_2)|=1$, say $c(V_1, V_2)=1$. If $|V_2|\geq 2$, then $1 \notin C(V_1)$ by Fact \ref{fa:3.1} (2).
Apply induction on $k$ within $V_1$, we have $|V(G)|=|V_1|+|V_2|\leq 2|V_1|\leq 2[gr_{k-1}(K_3 : F_{10})-1]< n$,
a contradiction. Thus $|V_2|=1$ and $|V_1|=n-1$. By Fact \ref{fa:3.1} (1), color 1 induces a subgraph $G^{(1)}$
such that each component is a $K_3$ or a star within $V_1$. In order to avoid a rainbow $K_3$ and
a monochromatic $P_4$ in color 1, there is only a single color between each pair of components.
By Lemma \ref{le:K3}, we have
$$gr_{k-1}(K_3 : K_3)=
\left\{
   \begin{aligned}
    &2\cdot 5^{(k-2)/2}+1,\ \ \ & \ & \ \ \mbox{if $k$ is even},\\
    &5^{(k-1)/2}+1,\ \ \ & \ & \ \ \mbox{if $k$ is odd}.
   \end{aligned}
   \right.$$
Since $|V_1|$ cannot be divided by three and a monochromatic $K_{1,3,3}$ contains a $F_{10}$,
$G^{(1)}$ contains at most $n_k=gr_{k-1}(K_3 : K_3)-2$ components each of which is a $K_3$.
For $0\leq i\leq n_k$, if there are exactly $n_k - i$ components each of which is a $K_3$, then there are at most
$i+1$ components each of which is a star to avoid a monochromatic $K_{2,2,3}$.
Thus after removing at most $2(n_k-i)+i+1$ vertices, there is no edge using color 1 within $V_1$.
Thus $|V(G)|=|V_1|+|V_2|\leq 2(n_k-i)+i+1+[gr_{k-1}(K_3 : F_{10})-1]+1< n$, a contradiction.

Therefore, we may assume $4 \leq m \leq 8$.
Let $r$ be the number of parts with at least three vertices, i.e.,
$|V_r|\geq 3$ and $|V_{r+1}|\leq 2$. Since $m\leq 8$ and $|V(G)|\geq 21$,
we have $1\leq r\leq m \leq 8$.
We divide the rest of the proof into four cases
based on the number $r$.

{\bf Case 1.} $r=1$.
\vspace{0.05cm}

For $i=1, 2$, let $A_i$ be the union of parts with color $i$ to $V_1$.
Since $m\geq 4$, we have $|A_1|\geq 1$, $|A_2|\geq 1$ and at least one of
$A_1$ and $A_2$ contains more that one vertices, say $|A_1|\geq 2$.
Then $1\notin C(V_1)$ by Fact \ref{fa:3.1} (2), and so $|V_1|\leq gr_{k-1}(K_3 : F_{10})-1$.
Since $2\geq |V_2|\geq \cdots \geq |V_m|$ and $m\leq 8$, $|A_1|+|A_2|\leq 14$.
If $|A_2|\geq 2$, then $2\notin C(V_1)$ and thus
$|V(G)|= |V_1|+|A_1|+|A_2|\leq gr_{k-2}(K_3 : F_{10})-1+14 <n$,
a contradiction. Thus $|A_2|=1$. Then $|A_1|\geq 3$ since otherwise
$|V(G)|= |V_1|+|A_1|+|A_2|\leq gr_{k-1}(K_3 : F_{10})-1+2+1 <n$.
By Fact \ref{fa:3.1} (2), $1\notin C(A_1)$, and thus $|A_1|\leq gr_{k-1}(K_3 : F_{10})-1$.
Then $|V(G)|= |V_1|+|A_1|+|A_2|\leq gr_{k-1}(K_3 : F_{10})-1+\mbox{min}\{gr_{k-1}(K_3 : F_{10})-1+1, 14\} <n$,
a contradiction.

{\bf Case 2.} $r=2$.
\vspace{0.05cm}

In this case, we may assume that $c(V_1, V_2)=1$. Then by Fact \ref{fa:3.1} (2), for $i=1, 2$,
$1\notin C(V_i)$ and thus $|V_i|\leq gr_{k-1}(K_3 : F_{10})-1$.
Since $m\leq 8$ and $2\geq |V_3|\geq \cdots \geq |V_m|$, if $|V_3 \cup \cdots \cup V_m|\geq 11$
then there is a monochromatic $K_{1, 2, 2}$, contradicting to the Fact \ref{fa:3.1} (3).
Hence, $|V_3 \cup \cdots \cup V_m|\leq 10$.
By Fact \ref{fa:3.1} (3), we also have that there is no vertex with color 1 to both
$V_1$ and $V_2$. Thus we can partition $V_3 \cup \cdots \cup V_m$ into three parts $A$, $B$
and $C$ such that, $c(A, V_1)=c(C, V_2)=1$ and $c(A, V_2)=c(C, V_1)=c(B, V_1\cup V_2)=2$.
If $|A\cup B|\geq 2$ and $|C\cup B|\geq 2$, then for $i=1, 2$, we have $2\notin C(V_i)$ and thus
$|V(G)|=|V_1|+|V_2|+|V_3 \cup \cdots \cup V_m|\leq 2[gr_{k-2}(K_3 : F_{10})-1]+10<n$,
a contradiction. Thus $|A\cup B|\leq 1$ or $|C\cup B|\leq 1$, say $|C\cup B|\leq 1$, so $|B|\leq 1$.

We first suppose $|B|=1$ and $|C|=0$.
Since $m\geq 4$, $|A|\geq 1$, and thus $2 \notin C(V_2)$. Since
$|A|=|V(G)|-|V_1|-|V_2|-|B|\geq n-[gr_{k-1}(K_3 : F_{10})-1]-[gr_{k-2}(K_3 : F_{10})-1]-1>2$,
we have $1, 2 \notin C(A)$.
Thus $|V(G)|=|V_1|+|V_2|+|A|+|B|+|C|\leq [gr_{k-1}(K_3 : F_{10})-1]+2[gr_{k-2}(K_3 : F_{10})-1]+1 < n$,
a contradiction.
Next we suppose $|B|=0$.
Since $m\geq 4$, we have $|A|\geq 1$ and $|C|= 1$.
If $|A|\geq 2$, then $2\notin C(V_2)$ and $|V_2|\leq gr_{k-2}(K_3 : F_{10})-1$.
Since $|A|=|V(G)|-|V_1|-|V_2|-|C|\geq n-[gr_{k-1}(K_3 : F_{10})-1]-[gr_{k-2}(K_3 : F_{10})-1]-1>2$,
we have $1, 2 \notin C(A)$, and thus $|V(G)|\leq [gr_{k-1}(K_3 : F_{10})-1]+2[gr_{k-2}(K_3 : F_{10})-1]+1 < n$,
a contradiction. Thus, $|A|=1$. By the minimality of $m$, we have $c(A, C)=2$.
By Fact \ref{fa:3.1} (1), color 2 induces a subgraph $G^{(2)}$
such that each component is a $K_3$ or a star within $V_1$. In order to avoid a rainbow $K_3$ and
a monochromatic $P_4$ in color 2, there is only a single color between each pair of components.
Since $1\notin C(V_1)$ and by Fact \ref{fa:3.1} (3), $G^{(2)}$ contains at most $n'_k=gr_{k-2}(K_3 : K_3)-1$
components each of which is a $K_3$ for avoiding a monochromatic $K_{3,3,3}$ that contains a $F_{10}$.
For $0\leq j\leq n'_k$, if there are exactly $n'_k - j$ components each of which is a $K_3$, then there are at most
$j$ components each of which is a star to avoid a monochromatic $K_{2,3,3}$.
Thus after removing at most $2(n'_k-j)+j$ vertices, there is no edge using color 2 within $V_1$.
Thus $|V(G)|=|V_1|+|V_2|+|A|+|C|\leq 2|V_1|+2\leq 2[2(n'_k-j)+j+(gr_{k-2}(K_3 : F_{10})-1)]+2< n$, a contradiction.

{\bf Case 3.} $r=3$.
\vspace{0.05cm}

By Fact \ref{fa:3.1} (3),, we may assume that $c(V_1, V_2\cup V_3)=1$
and $c(V_2, V_3)=2$. Thus $1 \notin C(V_1)\cup C(V_2)\cup C(V_3)$ and $2 \notin C(V_2)\cup C(V_3)$.
By inductive assumption, $|V_1|\leq gr_{k-1}(K_3 : F_{10})-1$, $|V_2|\leq gr_{k-2}(K_3 : F_{10})-1$
and $|V_3|\leq gr_{k-2}(K_3 : F_{10})-1$. Since $m\leq 8$ and $2\geq |V_4|\geq \cdots \geq |V_m|$,
we have $|V_4|\geq \cdots \geq |V_m|\leq 10$.
If $c(v, V_1)=1$ for some $v \in V_4 \cup \cdots \cup V_m$, then $c(v, V_2)\neq 1$
and $c(v, V_3)\neq 1$ by Fact \ref{fa:3.1} (3). But then $c(v, V_2)=c(v, V_3)=2$,
a contradiction. Thus $c(V_1, V_4 \cup \cdots \cup V_m)=2$.
Note that $|V_4 \cup \cdots \cup V_m|=|V(G)|-|V_1|-|V_2|-|V_3|\geq 5$.
By Fact \ref{fa:3.1} (2), $2\notin C(V_1)$ and $2\notin C(V_4 \cup \cdots \cup V_m)$.
Thus $|V(G)|\leq 3[gr_{k-2}(K_3 : F_{10})-1]+\mbox{min}\{gr_{k-1}(K_3 : F_{10})-1, 10\} < n$,
a contradiction.

{\bf Case 4.} $r\geq 4$.
\vspace{0.05cm}

By Fact \ref{fa:3.1} (3), there is no monochromatic $K_3$ in the subgraph $H$
of the reduced graph induced on the first four parts.
Thus $H$ is one of the two 2-coloring of $K_4$ with no monochromatic $K_3$.
By Fact \ref{fa:3.1} (2), $1, 2 \notin C(V_1)\cup C(V_2)\cup C(V_3)\cup C(V_4)$.
Thus for $i \in [4]$, $|V_i|\leq gr_{k-2}(K_3 : F_{10})-1$,
and so $|V_5 \cup \cdots \cup V_m|\geq 5$.
If $6 \leq m \leq 8$, then there is a monochromatic $K_{1, 2, 3}$, a contradiction.
Thus $m=5$ and $|V_5|\geq 5$.
In order two avoid a monochromatic $F_{10}$, $G$ is a blow-up of the unique 2-coloring
of $K_5$ with no monochromatic $K_3$.
By Fact \ref{fa:3.1} (2), $1, 2 \notin C(V_1)\cup C(V_2)\cup \cdots \cup C(V_5)$.
Apply induction on $k$, we have $|V(G)|\leq 5[gr_{k-2}(K_3 : F_{10})-1]< n$,
a contradiction.
\hfill$\blacksquare$

\section{Ramsey numbers and Gallai-Ramsey numbers for a class of unicyclic graph}
\label{sec:ch-3}

For $n\geq 3$, let $F_{2,n}$ denote the graph obtained by adding $n-2$ pendent edges
to a single vertex of $C_4$. The vertex with degree $n$ is called a center of $F_{2,n}$.
It is easy to see that $|V(F_{2,n})|=|E(F_{2,n})|=n+2$. Note that $F_{2,3}=F_{11}$.
In the following, we first give the 2-colored Ramsey number of $F_{2,n}$.

\noindent\begin{theorem}\label{th:3.1} For all integers $k\geq 1$ and $n\geq 3$,
$$r_2(F_{2,n})=
\left\{
   \begin{aligned}
    &2n-1,\ \ \  \ \ \ \ \ & \ & \ \ \mbox{if $n$ is even},\\
    &2n,\ \ \ & \ & \ \ \mbox{if $n$ is odd}.
   \end{aligned}
   \right.$$
\end{theorem}

\noindent {\bf Proof.} Let $n'=2n-\varepsilon$, where $\varepsilon=1$ if $n$ is even, and
$\varepsilon=0$ otherwise. Since $K_{1,n}$ is a subgraph of $F_{2,n}$
and $r_2(K_{1,n})=n'$ (proven in \cite{BuRo}), we have $r_2(F_{2,n})\geq n'$.

For the upper bound, suppose $G$ is a 2-coloring of $K_{n'}$ containing no
monochromatic $F_{2,n}$. Since $r_2(F_{2,3})=6$ (see \cite{Hen}), we may assume that $n\geq 4$.
Then $G$ contains a monochromatic $K_{1,n}$, say
using vertex set $\{u, v_1, v_2, \ldots, v_{n}\}$ and $c(u, \{v_1, v_2, \ldots, v_{n}\})=1$.
Let $V(G)\setminus \{u, v_1, v_2, \ldots, v_{n}\} = \{w_1, w_2, \ldots, w_{n'-n-1}\}$.
For any $i \in [n'-n-1]$, there are at most one edge using color 1 in $E(w_i, \{v_1, v_2, \ldots, v_{n}\})$
in order to avoid a monochromatic $F_{2,n}$.
Note that $n\geq 4$ and $n'-n-1 \geq 2$. If $c(w_i, \{v_1, v_2, \ldots, v_{n}\})=2$ for some $i \in [n'-n-1]$,
then it is easy to find a monochromatic $F_{2,n}$ in color 2, a contradiction.
Thus we may consider the following two cases based on $C(\{w_1, w_2\}, \{v_1, v_2, \ldots, v_{n}\})$.

Firstly, suppose that $c(\{w_1, w_2\}, \{v_1, v_2, \ldots, v_{n-1}\})=2$ and $c(\{w_1, w_2\}, v_{n})=1$.
In order to avoid a monochromatic $F_{2,n}$ in color 2, we have $c(\{w_1, w_2\}, u)=1$.
Then $uw_1v_nw_2u$ is a $C_4$ in color 1, which together with $\{v_1, v_2, \ldots, v_{n-2}\}$
forms a monochromatic $F_{2,n}$ centered at $u$, a contradiction.
Secondly, suppose that $c(\{w_1, w_2\}, \{v_1, v_2, \ldots, v_{n-2}\})=c(w_1v_n)=c(w_2v_{n-1})=2$
and $c(w_1v_{n-1})=c(w_2v_{n})=1$. Note that we also have $c(\{w_1, w_2\}, u)=1$ similarly as above.
If $c(v_1v_{n-1})=1$, then $uw_1v_{n-1}v_1u$ is a $C_4$ in color 1, which together
with $\{v_2, v_3, \ldots, v_{n-2}, v_{n}\}$ forms a monochromatic $F_{2,n}$ centered at $u$, a contradiction.
Thus $c(v_1v_{n-1})=2$, and by symmetry $c(\{v_1, v_2, \ldots, v_{n-2}\}, \{v_{n-1}, v_{n}\})=2$.
If $n=4$, then there is a $F_{2,4}$ in color 2 centered at $v_1$.
If $n\geq 5$, then $n'-n-1 \geq 4$ and we may assume $c(w_3v_{n-1})=2$ without loss
of generality. Then $v_{n-1}v_1v_nv_2v_{n-1}$ is a $C_4$ in color 2, which together with
$\{v_3, v_4, \ldots, v_{n-2}, w_2, w_3\}$ forms a monochromatic $F_{2,n}$ centered at $v_{n-1}$, a contradiction.
\hfill$\blacksquare$

\noindent\begin{theorem}\label{th:3.2}
~

(1) For all integers $k\geq 1$ and $n \in \{3,4\}$, $gr_k(K_3 : F_{2,n})=r_2(F_{2,n})+k-2$;

(2) For any integer $k\geq 3$, $gr_k(K_3 : F_{2,5})=k+9$;

(3) For all integers $k\geq 3$ and $n \geq 6$,
$$k(n-1)+2 \geq gr_k(K_3 : F_{2,n})\geq
\left\{
   \begin{aligned}
    &\frac{5n}{2}+k-6,\ \ \  \ \ \ \ \ & \ & \ \ \mbox{if $n$ is even},\\
    &\frac{5n-1}{2}+k-4,\ \ \ & \ & \ \ \mbox{if $n$ is odd}.
   \end{aligned}
   \right.$$
\end{theorem}

\noindent {\bf Proof.}
We first prove the lower bounds by constructing some appropriate Gallai-$k$-colorings
of complete graphs without monochromatic $F_{2,n}$.
For $n\in \{3, 4\}$, let $G_2$ be a 2-coloring of complete graph on $r_2(F_{2,n})-1$ vertices
using colors 1 and 2
without monochromatic $F_{2,n}$. Let $G'$ be the unique 2-coloring
of $K_5$ using colors 2 and 3 with no monochromatic $K_3$.
For $n\geq 5$ and $n$ even, let $G_3$ be the 3-coloring obtained by
substituting $K_{n/2}$ (using only color 1) into one vertex of $G'$ and
$K_{n/2-1}$ (using only color 1) into the other four vertices of $G'$.
For $n\geq 5$ and $n$ odd, let $G_3$ be the 3-coloring obtained by
substituting $K_{(n-1)/2}$ (using only color 1) into five vertex of $G'$.
Suppose $i < k$ and we have constructed $G_i$ containing no rainbow $K_3$
and no monochromatic $F_{2,n}$. We construct $G_{i+1}$ by adding a new vertex to
$G_{i}$ such that all the new edges are colored by color $i+1$.
Finally, we obtain a $k$-coloring $G_{k}$ containing no rainbow $K_3$
and no monochromatic $F_{2,n}$. For $n \in \{3,4\}$, $|V(G_k)|=r_2(F_{2,n})+k-3$.
For $n\geq 6$, we have
$$|V(G_k)|=
\left\{
   \begin{aligned}
    &\frac{5n}{2}+k-7,\ \ \  \ \ \ \ \ & \ & \ \ \mbox{if $n$ is even},\\
    &\frac{5n-1}{2}+k-5,\ \ \ & \ & \ \ \mbox{if $n$ is odd}.
   \end{aligned}
   \right.$$
For $n=5$, we construct $G_{k+1}$ by adding a new vertex to
$G_{k}$ such that all the new edges are colored by color $1$.
Then $G_{k+1}$ contains no rainbow $K_3$ and no monochromatic $F_{2,5}$, and
$|V(G_k)|=k+8$.

For the upper bound, let
$$n_{k}'=
\left\{
   \begin{aligned}
    &r_2(F_{2,n})+k-2,\ \ \  \ \ \ \ \ & \ & \ \ \mbox{if $n \in \{3,4\}$},\\
    &k+9,\ \ \  \ \ \ \ \ & \ & \ \ \mbox{if $n=5$},\\
    &k(n-1)+2,\ \ \ & \ & \ \ \mbox{if $n \geq 6$}.
   \end{aligned}
   \right.$$
We will prove that every Gallai-$k$-coloring of $K_{n_{k}'}$ contains a monochromatic $F_{2,n}$
for all integers $k\geq 1$ and $n\geq 3$.
The case $k=1$ is trivial and the case $k=2$ holds by Theorem \ref{th:3.1},
so we may assume $k\geq 3$ in the following.
For a contradiction, suppose $G$ is a $k$-coloring of $K_{n'}$ containing no
rainbow $K_3$ and no monochromatic $F_{2,n}$.
We choose such a $G$ with $k$ minimal.
Let $V_1, V_2, \ldots, V_m$ be a Gallai-partition with $|V_1|\geq |V_2| \geq \cdots \geq |V_m|$
and $2\leq m\leq r_2(F_{2,n})-1$. Suppose colors 1 and 2 are the two colors used between the parts.

We first suppose that $m\geq 4$. In this case, we have $|V(G)\setminus V_1|\geq 3$.
In order to avoid a monochromatic $F_{2,n}$ and since $k\geq 3$, we have $2\leq |V_{1}|\leq n-1$.
For $i=1, 2$, let $A_i$ be the union of parts with color $i$ to $V_1$. If $n\notin \{4, 5\}$, then
$\mbox{max}\{|A_1|, |A_2|\}\geq \lceil\frac{|V(G)|-|V_1|}{2}\rceil \geq \lceil\frac{n_{k}'-(n-1)}{2}\rceil \geq n$,
which implies that there is a monochromatic $F_{2,n}$, a contradiction.
If $n=5$, we must have $|V_1|=|A_1|=|A_2|=4$ and $k=3$ by a similar argument.
For any edge $uv \in E(A_1, A_2)$, we have $c(uv)\in \{1, 2\}$, resulting in a
monochromatic $F_{2,5}$, a contradiction.
If $n=4$, we have $2\leq |V_1|\leq 3$ and
$3\geq \mbox{max}\{|A_1|, |A_2|\}\geq \lceil\frac{k+5-3}{2}\rceil \geq 3$, say $|A_1|=3$.
If $|V_1|=2$, then $|A_2|=3$ and $k=3$. Let $A_{1}=\{u_1, u_2, u_3\}$
and $A_2=\{v_1, v_2, v_3\}$. In order to avoid a monochromatic $F_{2,4}$,
there are at most one edge using color 1 between $u_i$ and $A_2$ for each $i\in [3]$.
Thus there are at most three edges using color 1 between $A_1$ and $A_2$.
Then there are at least six edges using color 2 between $A_1$ and $A_2$,
and in particular, there exist a vertex $v_i\in A_{2}$ such that there are at least two edges
using color 2 between $v_i$ and $A_1$, resulting in a monochromatic $F_{2,4}$, a contradiction.
If $|V_1|=3$, then $2\leq |A_2|\leq 3$.
For avoiding a monochromatic $F_{2,4}$ in color 1, $1\notin C(A_1, A_2)$, i.e.,
$c(A_1, A_2)=2$, resulting in a monochromatic $F_{2,4}$ in color 2, a contradiction.

Therefore, we have $2\leq m\leq 3$. In this case, there exists a bipartition $V_1$,
$V_2$ of $V(G)$ with $|C(V_1, V_2)|=1$, say $c(V_1, V_2)=1$.
Since $k\geq 3$, $|V_1|\geq \lceil \frac{|V(G)|}{2}\rceil \geq n$.
Thus $|V_2|=1$ for avoiding a monochromatic $F_{2,n}$ in color 1.
Moreover, if $1\notin C(V_1)$, then $|V(G)|=|V_1|+|V_2|\leq (n_{k-1}'-1)+1< n_{k}'$
by the minimality of $k$, a contradiction. Thus $1 \in C(V_1)$, say $v_1, v_2 \in V_1$
and $c(v_1v_2)=1$. In order to avoid a monochromatic $F_{2,n}$ in color 1,
$1\notin C(\{v_1, v_2\}, V_1\setminus \{v_1, v_2\})$.
Since $G$ is rainbow $K_3$-free, we have $c(v_1v)=c(v_2v)$ for any $v\in V_1\setminus \{v_1, v_2\}$.
For avoiding a monochromatic $F_{2,n}$, there are at most $n-1$ vertices in $V_1\setminus \{v_1, v_2\}$
with a single color to $\{v_1, v_2\}$.
Note that $|V_1\setminus \{v_1, v_2\}|=n_{k}'-3$ and $|C(\{v_1, v_2\}, V_1\setminus \{v_1, v_2\})|=k-1$.

If $n\geq 6$, then there are at least $\lceil \frac{n_{k}'-3}{k-1}\rceil \geq n$ vertices in $V_1\setminus \{v_1, v_2\}$
with a single color to $\{v_1, v_2\}$, i.e., there is a monochromatic $F_{2,n}$, a contradiction.

If $n=3$, then there exist four vertices $u_1, u_2, u_3, u_4 \in V_1\setminus \{v_1, v_2\}$
such that $c(\{v_1, v_2\},$ $\{u_1, u_2\})=i$ and $c(\{v_1, v_2\}, \{u_3, u_4\})=j$, where
$2\leq i< j\leq k$. Since $G$ is rainbow $K_3$-free, we have
$C(\{u_1, u_2\}, \{u_3, u_4\})\subseteq \{i, j\}$. Then there is a monochromatic $F_{2,3}$
in color $i$ or $j$, a contradiction.

If $n=4$, then we may consider the following two cases.
First, there exist five vertices $u_1, u_2, \ldots, u_5 \in V_1\setminus \{v_1, v_2\}$
such that $c(\{v_1, v_2\},$ $\{u_1, u_2, u_3\})=i$ and $c(\{v_1, v_2\}, \{u_4, u_5\})$ $=j$, where
$2\leq i< j\leq k$. Since $G$ is rainbow $K_3$-free, we have
$C(\{u_1, u_2, u_3\}, \{u_4, u_5\})\subseteq \{i, j\}$.
In order to avoid a monochromatic $F_{2,4}$ in color $j$, there are at least two edges using color $i$
between $u_4$ (resp., $u_5$) and $\{u_1, u_2, u_3\}$. Thus there exists a vertex, say $u_1$, such that
$c(u_1, \{u_4, u_5\})=i$, resulting in a monochromatic $F_{2,4}$,
a contradiction. Second, there exist six vertices $u_1, u_2, \ldots, u_6 \in V_1\setminus \{v_1, v_2\}$
such that $c(\{v_1, v_2\},$ $\{u_1, u_2\})=i$, $c(\{v_1, v_2\}, \{u_3, u_4\})=j$ and
$c(\{v_1, v_2\}, \{u_5, u_6\})=l$, where $2\leq i< j < l \leq k$. Without loss of
generality,we may assume that $c(u_1u_3)=i$.
Then $c(u_1u_4)=j$ and $c(u_1, \{u_5, u_6\})=l$.
Since $c(u_4u_5)\in \{j, l\}$, there is a monochromatic $F_{2,4}$ in color $j$ or $l$,
a contradiction.

If $n=5$, then by the pigeonhole principle, we may consider the following two
cases. First, there exist nine vertices $u_1, u_2, \ldots, u_9 \in V_1\setminus \{v_1, v_2\}$
such that $c(\{v_1, v_2\},$ $\{u_1, u_2, u_3, u_4\})=i$, $c(\{v_1, v_2\}, \{u_5, u_6, u_7\})=j$ and
$c(\{v_1, v_2\}, \{u_8, u_9\})=l$, where $2\leq i< j < l \leq k$.
In order to avoid a monochromatic $F_{2,5}$ in color $j$ (resp., $l$), there are at most two edges using color $j$
(resp., $l$) between each $u\in \{u_5, u_6, u_7\}$ (resp., $\{u_8, u_9\}$) and $\{u_1, u_2, u_3, u_4\}$.
Thus, there are at least ten edges using color $i$ between $\{u_1, u_2, u_3, u_4\}$ and $\{u_5, u_6, \ldots, u_9\}$,
and in particular, there exist a vertex $u\in \{u_1, u_2, u_3, u_4\}$ such that there are at least three edges
using color $i$ between $u$ and $\{u_5, u_6, \ldots, u_9\}$, resulting in a monochromatic $F_{2,5}$, a contradiction.
Second, there exist eight vertices $u_1, u_2, \ldots, u_8 \in V_1\setminus \{v_1, v_2\}$
such that $c(\{v_1, v_2\},$ $\{u_1, u_2\})=i$, $c(\{v_1, v_2\},$ $\{u_3, u_4\})=j$, $c(\{v_1, v_2\}, \{u_5, u_6\})=s$ and
$c(\{v_1, v_2\}, \{u_7, u_8\})=t$, where $2\leq i< j < s< t \leq k$.
Since there are at most two edges using color $i$ in $E(u_1, \{u_3, u_4, \ldots, u_8\})$,
we may assume that $c(u_1, \{u_7, u_8\})=t$ without loss of generality.
Then there is at most one edge using color $t$ in $E(u_7, \{u_3, u_4, u_5, u_6\})$,
we may assume that $c(u_7, \{u_5, u_6\})=s$. For avoiding a rainbow $K_3$,
we have $c(u_1, \{u_5, u_6\})=s$. For avoiding a monochromatic $F_{2,5}$ in color $s$, we
have $c(u_8, \{u_5, u_6\})=t$, resulting in a monochromatic $F_{2,5}$ in color $t$ with center $u_8$,
a contradiction.
\hfill$\blacksquare$

\section{Proof of Theorem \ref{th:main3}}
\label{sec:ch-4}

We begin with the following Ramsey numbers which will be used in the proofs.

\noindent\begin{lemma}\label{le:4.1} {\normalfont (\cite{Hen})}
$r_2(F_{12})=r_2(F_{13})=10$.
\end{lemma}

Theorem \ref{th:main3} follows from the following two lemmas immediately.

\noindent\begin{lemma}\label{le:4.2}
For any integer $k\geq 1$, $H\in \{F_{12}, F_{13}\}$,
$$gr_k(K_3 : H)>
\left\{
   \begin{aligned}
    &9\cdot 5^{(k-2)/2},\ \ \ \ \ & \ & \ \mbox{if $k$ is even},\\
    &4\cdot 5^{(k-1)/2},\ \ \ \ \ & \ & \ \mbox{if $k$ is odd}.
   \end{aligned}
   \right.$$
\end{lemma}

\noindent {\bf Proof.} If $k$ is even, then let $G_2$ be a 2-colored $K_9$
containing no monochromatic $H$ using colors 1 and 2. Suppose that $2i < k$ and we have constructed
a $2i$-coloring $G_{2i}$ of $K_{n_{2i}}$ containing neither rainbow $K_3$ nor monochromatic $H$,
where $n_{2i}=9\cdot 5^{(2i-2)/2}$. Let $G'$ be a 2-colored $K_5$ using colors $2i+1$ and $2i+2$
which contains no monochromatic $K_3$, i.e., colors $2i+1$ and $2i+2$ induce two monochromatic
$C_5$. We construct $G_{2i+2}$ by substituting five copies of $G_{2i}$ into vertices of $G'$,
i.e., $G_{2i+2}$ is a blow-up of $G'$. Finally, we obtain a $k$-coloring $G_k$ of $K_n$
containing neither rainbow $K_3$ nor monochromatic $H$, where $n=9\cdot 5^{(k-2)/2}$.

If $k$ is odd, then let $G_1$ be a monochromatic $K_4$ using color 1.
Suppose that $2i-1 < k$ and we have constructed
a $(2i-1)$-coloring $G_{2i-1}$ of $K_{n_{2i-1}}$ containing neither rainbow $K_3$ nor monochromatic $H$,
where $n_{2i-1}=4\cdot 5^{(2i-2)/2}$. Let $G''$ be a 2-colored $K_5$ using colors $2i$ and $2i+1$
which contains no monochromatic $K_3$. We construct $G_{2i+1}$ by
substituting five copies of $G_{2i-1}$ into vertices of $G''$.
Finally, we obtain a $k$-coloring $G_k$ of $K_n$ containing neither
rainbow $K_3$ nor monochromatic $H$, where $n=4\cdot 5^{(k-1)/2}$.
\hfill$\blacksquare$

\noindent\begin{lemma}\label{le:4.3}
For any integer $k\geq 1$, $H\in \{F_{12}, F_{13}\}$,
$$gr_k(K_3 : H)\leq
\left\{
   \begin{aligned}
    &9\cdot 5^{(k-2)/2}+1,\ \ \  & \ & \ \ \mbox{if $k$ is even},\\
    &4\cdot 5^{(k-1)/2}+1,\ \ \ & \ & \ \ \mbox{if $k$ is odd}.
   \end{aligned}
   \right.$$
\end{lemma}

\noindent {\bf Proof.} We prove the statement by induction on $k$.
The case $k=1$ is trivial and the case $k=2$ holds by Lemma \ref{le:4.1},
so we may assume $k \geq 3$. Suppose $G$ is a Gallai-$k$-coloring of $K_{n}$
containing no monochromatic $H$, where
$$n=
\left\{
   \begin{aligned}
    &9\cdot 5^{(k-2)/2}+1,\ \ \  & \ & \ \ \mbox{if $k$ is even},\\
    &4\cdot 5^{(k-1)/2}+1,\ \ \ & \ & \ \ \mbox{if $k$ is odd}.
   \end{aligned}
   \right.$$
Let $V_1, V_2, \ldots, V_m$ be a Gallai-partition with $|V_1|\geq |V_2| \geq \cdots \geq |V_m|$
and $2 \leq m \leq 9$. We choose such a partition with $m$ minimum.
Suppose colors 1 and 2 are the two colors used between the parts.
We first have the following simple facts since $G$ is monochromatic $H$-free.

\noindent\begin{fact}\label{fa:4.1} For any two non-empty parts $V_i$, $V_j$ with $c(V_i, V_j)=\alpha$, we have

(1) there is no monochromatic $P_4$ using color $\alpha$ within $V_i$ and $V_j$, respectively;

(2) if $|V_i|\geq 3$ and $|V_j|\geq 3$, then $\alpha \notin C(V_i)$ and $\alpha \notin C(V_i)$, respectively;

(3) if $|V_i|\geq 3$ and $|V_j|\geq 2$, then there is no monochromatic $P_3$ using color $\alpha$ within $V_i$;

(4) if $|V_i|\geq 4$, then there is no monochromatic $K_3$ using color $\alpha$ within $V_i$;

(5) if $|V_i|\geq 2$ and $|V_j|\geq 2$, then there is no vertex with color $\alpha$ to both
$V_i$ and $V_j$ in $V(G)\setminus (V_i \cup V_j)$.
\end{fact}

\noindent\begin{fact}\label{fa:4.2} For any two non-empty parts $V_i$, $V_j$ with $c(V_i, V_j)=\alpha$, we have

(1) if $H=F_{12}$ and $|V_i|\geq 4$, then there is no monochromatic $P_3$ using color $\alpha$ within $V_i$;

(2) if $H=F_{13}$, $|V_i|\geq 3$ and $|V_j|\geq 2$, then $\alpha \notin C(V_i)$.
\end{fact}

We now claim that $4 \leq m \leq 9$. In fact, if $2\leq m \leq 3$, then there exists a bipartition
$V_1$ and $V_2$ of $V(G)$ with $|C(V_1, V_2)|=1$, say $c(V_1, V_2)=1$. If $|V_2|\geq 3$, then $1 \notin C(V_1)$ by Fact \ref{fa:4.1} (2).
Apply induction on $k$ within $V_1$, we have $|V(G)|=|V_1|+|V_2|\leq 2|V_1|\leq 2[gr_{k-1}(K_3 : H)-1]< n$,
a contradiction. If $|V_2|= 2$, then by Fact \ref{fa:4.1} (3) color 1 induces a matching within $V_1$.
If $|V_2|= 1$, then by Fact \ref{fa:4.1} (1) and (4), color 1 induces a subgraph such that each component
is a star within $V_1$. Thus when $1\leq |V_2|\leq 2$, we can partition $V_1$ into two subsets $V_1'$
and $V_1''$ such that there is no color 1 within $V_1'$ and $V_1''$, respectively.
By inductive assumption, $|V(G)|=|V_1'|+|V_1''|+|V_2|\leq 2[gr_{k-1}(K_3 : H)-1]+2< n$, a contradiction.

Therefore, we may assume $4 \leq m \leq 9$.
Let $r$ be the number of parts with at least three vertices, i.e.,
$|V_r|\geq 3$ and $|V_{r+1}|\leq 2$. Since $m\leq 9$ and $|V(G)|\geq 21$,
we have $1\leq r\leq m \leq 9$.
We divide the rest of the proof into four cases
based on the number $r$.

{\bf Case 1.} $r\geq 4$.
\vspace{0.05cm}

By Fact \ref{fa:4.1} (5), there is no monochromatic $K_3$ in the subgraph $R$
of the reduced graph induced on the first four parts.
Thus $R$ is one of the two 2-coloring of $K_4$ with no monochromatic $K_3$.
If $m=5$, then the reduced graph is the unique 2-coloring of $K_5$ with no monochromatic $K_3$.
Thus for $4\leq m\leq 5$, we have $1, 2\notin C(V_1)$ by Fact \ref{fa:4.1} (2),
and so $|V(G)|\leq 5|V_1|\leq 5[gr_{k-2}(K_3 : H)-1]< n$, a contradiction.
If $6\leq m\leq 9$, then $1\geq |V_5|\geq \cdots \geq |V_m|$ otherwise there is
a monochromatic $K_{1,2,3}$.
Since $n=|V(G)|=\sum^{4}_{i=1}|V_i| + \sum^{m}_{j=5}|V_j| \leq 4|V_1|+5 \leq 4[gr_{k-2}(K_3 : H)-1]+5$,
we have $n=21$, $k=3$, $|V_i|=4$ for $i\in [4]$ and $|V_j|=1$ for $j\in \{5, 6, \ldots, 9\}$.
Note that there is a monochromatic $K_3$ in subgraph of the reduced graph induced on $V_1, V_2, \ldots, V_6$,
i.e., there is a monochromatic $K_{1,3,3}$ or $K_{1,1,3}$ in $G$.
If the former holds, then we are done.
If the latter holds, then there is a monochromatic $F_{12}$ using six edges of the $K_{1,3,3}$, and a
monochromatic $F_{13}$ using five edges of the $K_{1,3,3}$ and one edge between the first four parts, a contradiction.

{\bf Case 2.} $r=3$.
\vspace{0.05cm}

In order to avoid a monochromatic $K_{3,3,3}$, we may assume that $c(V_1, V_2\cup V_3)=1$
and $c(V_2, V_3)=2$. Then $1 \notin C(V_1)\cup C(V_2)\cup C(V_3)$ and $2 \notin C(V_2)\cup C(V_3)$.
By inductive assumption, $|V_1|\leq gr_{k-1}(K_3 : H)-1$, $|V_2|\leq gr_{k-2}(K_3 : H)-1$
and $|V_3|\leq gr_{k-2}(K_3 : H)-1$, and so $|V_4 \cup \cdots \cup V_m|\geq 4$.
Note that $|V_4 \cup \cdots \cup V_m|\leq 10$, otherwise
there is a monochromatic $K_{1,2,2}$ since $m\leq 9$ and $2\geq |V_4|\geq \cdots \geq |V_m|$.
If $c(v, V_1)=1$ for some $v \in V_4 \cup \cdots \cup V_m$, then $c(v, V_2)\neq 1$
and $c(v, V_3)\neq 1$ by Fact \ref{fa:4.1} (5). But then $c(v, V_2)=c(v, V_3)=2$,
a contradiction. Thus $c(V_1, V_4 \cup \cdots \cup V_m)=2$.
Recall that we have $|V_4 \cup \cdots \cup V_m|\geq 4$,
and thus $2\notin C(V_1)$ and $2 \notin C(V_4 \cup \cdots \cup V_m)$ by Fact \ref{fa:4.1} (2).
Hence, $|V_1|\leq gr_{k-2}(K_3 : H)-1$  and $|V_4 \cup \cdots \cup V_m|\leq \mbox{min}\{10, gr_{k-1}(K_3 : H)-1\}$.
Since $n=|V(G)|=\sum^{3}_{i=1}|V_i| + \sum^{m}_{j=4}|V_j| \leq 3[gr_{k-2}(K_3 : H)-1]+\mbox{min}\{10, gr_{k-1}(K_3 : H)-1\}$,
we have $n=21$, $k=3$ and $|V_4 \cup \cdots \cup V_m|=9$.
For any vertex $v \in V_4 \cup \cdots \cup V_m$, at least one of $c(v, V_2)$
and $c(v, V_3)$ is 1. Thus there are at least five vertices in $V_4 \cup \cdots \cup V_m$,
say $\{x_1, x_2, \ldots, x_5\}$, such that $c(\{x_1, x_2, \ldots, x_5\}, V_2)=1$
or $c(\{x_1, x_2, \ldots, x_5\}, V_3)=1$. Then $1\notin \{x_1, x_2, \ldots, x_5\}$,
a contradiction.

{\bf Case 3.} $r=1$.
\vspace{0.05cm}

For $i=1, 2$, let $A_i$ be the union of parts with color $i$ to $V_1$.
Then $|A_1|+|A_2|\leq 16$.
If $|A_1|\geq 3$ and $|A_2|\geq 3$, then $1, 2 \notin C(V_1)$ by Fact \ref{fa:4.1} (2).
By inductive assumption, $|V(G)|= |V_1|+|A_1|+|A_2|\leq gr_{k-2}(K_3 : H)-1+16<n$, a contradiction.
If $|A_1|\geq 3$ and $|A_2|\leq 2$, then $1 \notin C(V_1)$ and $1 \notin C(A_1)$ by Fact \ref{fa:4.1} (2).
By inductive assumption, $|V(G)|= |V_1|+|A_1|+|A_2|\leq 2[gr_{k-1}(K_3 : H)-1]+2<n$, a contradiction.
Thus $|A_1|\leq 2$, and by symmetry $|A_2|\leq 2$.
Since $m\geq 4$, we have $|A_1|\geq 1$, $|A_2|\geq 1$ and $\mbox{max}\{|A_1|, |A_2|\}=2$, say $|A_1|= 2$.

If $H=F_{12}$, then there is no monochromatic $P_3$ using color 1 (resp., color 2) within $V_1$ by Fact
\ref{fa:4.2} (1). Thus there is no rainbow $P_3$ using colors 1 and 2, since otherwise there is a
monochromatic $P_3$ using color 1 or 2 since $G$ is rainbow $K_3$-free.
Hence, colors 1 and 2 induce a matching within $V_1$.
Then we can remove at most $\lfloor \frac{|V_1|}{2}\rfloor$ vertices such that there is no color 1
and no color 2 within $V_1$.
By inductive assumption, $|V(G)|= |V_1|+|A_1|+|A_2|\leq 2[gr_{k-2}(K_3 : F_{12})-1]+4<n$, a contradiction.

If $H=F_{13}$, then $1\notin C(V_1)$ by Fact \ref{fa:4.2} (2).
Now $|V(G)|= |V_1|+|A_1|+|A_2|\leq gr_{k-1}(K_3 : F_{13})-1+4<n$, a contradiction.

{\bf Case 4.} $r=2$.
\vspace{0.05cm}

In this case, we may assume that $c(V_1, V_2)=1$. Then by Fact \ref{fa:4.1} (2), for $i=1, 2$,
$1\notin C(V_i)$ and thus $|V_i|\leq gr_{k-1}(K_3 : H)-1$.
Since $m\leq 9$ and $2\geq |V_3|\geq \cdots \geq |V_m|$, if $|V_3 \cup \cdots \cup V_m|\geq 12$
then there is a monochromatic $K_{1, 2, 2}$, contradicting to the Fact \ref{fa:4.1} (5).
Hence, $|V_3 \cup \cdots \cup V_m|\leq 11$.
By Fact \ref{fa:4.1} (5), we also have that there is no vertex with color 1 to both
$V_1$ and $V_2$. Thus we can partition $V_3 \cup \cdots \cup V_m$ into three parts $A$, $B$
and $C$ such that, $c(A, V_1)=c(C, V_2)=1$ and $c(A, V_2)=c(C, V_1)=c(B, V_1\cup V_2)=2$.
Note that we have $|A\cup B|\leq 2$ or $|B\cup C|\leq 2$, otherwise $2\notin C(V_1)$ and $2\notin C(V_2)$
by Fact \ref{fa:4.1} (2), and then $|V(G)|= |V_1|+|V_2|+|A|+|B|+|C|\leq 2[gr_{k-2}(K_3 : H)-1]+11<n$, a contradiction.
Without loss of generality, let $|B\cup C|\leq 2$, so $|B|\leq 2$, $|C|\leq 2$.

If $H=F_{12}$, then
we first suppose $|B|=2$ and $|C|=0$.
Since $m\geq 4$, $|A|\neq 0$.
If $|A|\geq 3$, then $2\notin C(V_2)$ and $1, 2 \notin C(A)$.
Thus $|V(G)|= |V_1|+|V_2|+|A|+|B|+|C|\leq [gr_{k-1}(K_3 : F_{12})-1]+2[gr_{k-2}(K_3 : F_{12})-1]+2<n$, a contradiction.
Hence, $|A|\leq 2$, but then
$|V(G)|= |V_1|+|V_2|+|A|+|B|+|C|\leq [gr_{k-1}(K_3 : F_{12})-1]+[gr_{k-2}(K_3 : F_{12})-1]+4<n$, a contradiction.
Next, we suppose $|B|=1$ and $|C|\leq 1$.
If $|A|\geq 2$, then we can derive a contradiction similarly as above.
Thus $|A|\leq 1$.
If $|A|=0$ or $|C|=0$, then $|V(G)|\leq 2[gr_{k-1}(K_3 : F_{12})-1]+2<n$.
Hence, $|A|=|C|=1$.
By Fact \ref{fa:4.1} (3), color 2 induces a matching within $V_1$ and $V_2$.
For $i=1, 2$, there is no color 2 within $V_i$ after removing at most $\lfloor \frac{|V_i|}{2}\rfloor$ vertices.
Thus $|V(G)|\leq 4[gr_{k-2}(K_3 : F_{12})-1]+4<n$, a contradiction.
Finally, we suppose $|B|=0$ and $|C|\leq 2$.
If $|A|\geq 3$, then $2\notin C(V_2)$ and $1, 2 \notin C(A)$ but
then $|V(G)|\leq [gr_{k-1}(K_3 : F_{12})-1]+2[gr_{k-2}(K_3 : F_{12})-1]+2<n$.
Thus $|A|\leq 2$.
Note $|V_2|\geq 4$ otherwise $|V(G)|\leq [gr_{k-1}(K_3 : F_{12})-1]+3+4<n$.
Since $m\geq 4$, we have $|A|\neq 0$ and $|C|\neq 0$.
By Fact \ref{fa:4.2} (1), color 2 induces a matching within $V_1$ and $V_2$.
Thus $|V(G)|\leq 4[gr_{k-2}(K_3 : F_{12})-1]+4<n$, a contradiction.

If $H=F_{13}$, then $|A\cup B|\leq 1$ or $|B\cup C|\leq 1$, otherwise $2\notin C(V_1)$ and $2\notin C(V_2)$
by Fact \ref{fa:4.2} (2), and then $|V(G)|= |V_1|+|V_2|+|A|+|B|+|C|\leq 2[gr_{k-2}(K_3 : F_{13})-1]+11<n$, a contradiction.
Without loss of generality, let $|B\cup C|\leq 1$, so $|B|\leq 1$, $|C|\leq 1$.
We first suppose $|B|=1$ and $|C|=0$.
Since $m\geq 4$, $|A|\neq 0$. Then $2\notin C(V_2)$ by Fact \ref{fa:4.2} (2),
and thus $|A|=|V(G)|-|V_1|-|V_2|-|B|-|C|\geq n-[gr_{k-1}(K_3 : F_{13})-1]-[gr_{k-2}(K_3 : F_{13})-1]-1>3$.
Then $1, 2 \notin C(A)$, and so
$|V(G)|= |V_1|+|V_2|+|A|+|B|+|C|\leq [gr_{k-1}(K_3 : F_{13})-1]+2[gr_{k-2}(K_3 : F_{13})-1]+1<n$, a contradiction.
Therefore, $|B|=0$ and $|C|\leq 1$. Since $m\geq 4$, we have $|A|\geq 1$ and $|C|=1$.
If $|A|\geq 3$, then $2\notin C(V_2)$ and $1, 2 \notin C(A)$.
Thus $|V(G)|\leq [gr_{k-1}(K_3 : F_{13})-1]+2[gr_{k-2}(K_3 : F_{13})-1]+1<n$, a contradiction.
If $|A|= 2$, then $2\notin C(V_2)$.
Thus $|V(G)|\leq [gr_{k-1}(K_3 : F_{13})-1]+[gr_{k-2}(K_3 : F_{13})-1]+3<n$, a contradiction.
Hence $|A|=|C|=1$.
Then $|V(G)|\leq 2[gr_{k-1}(K_3 : F_{13})-1]+2<n$, a contradiction.
\hfill$\blacksquare$


\begin{thebibliography}{99}


\bibitem{BuRo} S.A. Burr and J.A. Roberts, On Ramsey numbers for stars,
   {\it Utilitas Mathematica} (1973), 217--220.

\bibitem{ChLP} M. Chen, Y.S. Li and C.P. Pei, Gallai-Ramsey numbers of odd cycles and complete bipartite graphs,
   {\it Graphs Combin.} (2018). https://doi.org/10.1007/s00373-018-1931-7

\bibitem{ChGr} F.R.K. Chung and R. Graham, Edge-colored complete graphs with precisely colored subgraphs,
{\it Combinatorica} {\bf 3} (1983), 315--324.

\bibitem{FGJM} R.J. Faudree, R. Gould, M. Jacobson and C. Magnant, Ramsey numbers in rainbow
   triangle free colorings, {\it  Australas. J. Combin.} {\bf 46} (2010), 269--284.

\bibitem{FuMa1} S. Fujita and C. Magnant, Gallai-Ramsey numbers for cycles, {\it Discrete Math.} {\bf 311} (2011), 1247--1254.

\bibitem{FuMa2} S. Fujita and C. Magnant, Extensions of Gallai-Ramsey results, {\it J. Graph Theory} {\bf 70} (2012), 404--426.

\bibitem{FuMO1} S. Fujita, C. Magnant and K. Ozeki, Rainbow generalizations of Ramsey
   theory: a survey, {\it Graphs Combin.} {\bf 1} (2010), 1--30.

\bibitem{FuMO2} S. Fujita, C. Magnant and K. Ozeki, Rainbow generalizations of Ramsey
   theory: a dynamic survey, {\it Theo. Appl. Graphs} {\bf 0} (2014), Article 1.

\bibitem{Gallai} T. Gallai, Transitiv orientierbare Graphen, {\it Acta Math. Acad. Sci. Hungar} {\bf 18} (1967), 25--66.

\bibitem{GSSS} A. Gy\'{a}rf\'{a}s, G. N. S\'{a}rk\"{o}zy, A. Seb\H{o} and S. Selkow, Ramsey-type results for
   Gallai colorings, {\it J. Graph Theory} {\bf 64} (2010), 233--243.

\bibitem{GySi} A. Gy\'{a}rf\'{a}s and G. Simonyi, Edge colorings of complete graphs without
   tricolored triangles, {\it J. Graph Theory} {\bf 46} (2004), 211--216.

\bibitem{HMOT} M. Hall, C. Magnant, K. Ozeki and M. Tsugaki, Improved upper bounds
   for Gallai--Ramsey numbers of paths and cycles, {\it J. Graph Theory} {\bf 75} (2014), 59--74.

\bibitem{Hen} G.R.T. Hendry, Ramsey numbers for graphs with five vertices, {\it J. Graph Theory} {\bf 13} (1989), 245--248.

\bibitem{MaPr} F. Maffray and M. Preissmann, A translation of Tibor Gallai's paper: transitiv orientierbare Graphen, In: {\it J.L. Ramirez-Alfonsin, B.A.
   Reed (Eds.), Perfect Graphs, Wiley, New York} (2001), pp. 25--66.

\bibitem{Mag}  C. Magnant, Personal communication.

\bibitem{SoZh} Z.-X. Song and J.M. Zhang, A conjecture on Gallai-Ramsey numbers of even cycles and paths, arXiv:1803.07963 (2018).

\bibitem{WMMZ} Z. Wang, Y.P. Mao, C. Magnant and J.Y. Zou, Ramsey and Gallai-Ramsey numbers for two classes
   of unicyclic graphs, arXiv:1809.10298 (2018).

\bibitem{WMNX} H.B. Wu, C. Magnant, P.S. Nowbandegani and S.M. Xia, All partitions have small parts--Gallai-Ramsey numbers of bipartite graphs,
   {\it Discrete Appl. Math.} in press (2018).

\end{thebibliography}
\end{document}